\documentclass[12pt]{amsart}
\usepackage{latexsym,amssymb,amsfonts,amsmath}

\begin{document}

\begin{center}
{\Large \bf On Finite Groups of \vspace*{0.5em} \\
Symmetries of Surfaces} 
\end{center}

\vspace*{.4in}

\begin{center}
{\large\bf J\"{u}rgen M\"{u}ller and Siddhartha Sarkar}
\end{center}

\vspace*{.4in}

\hrule \smallskip
{\small \noindent{\bf Abstract.}
The genus spectrum of a finite group $G$ is the set of all $g\geq 2$ 
such that $G$ acts faithfully and orientation-preserving on a closed 
compact orientable surface of genus $g$. This article is an overview 
of some results relating the genus spectrum of $G$ to its group theoretical
properties.
In particular, the arithmetical properties of genus spectra
are discussed, and explicit results are given on the $2$-groups 
of maximal class, certain sporadic simple groups and
a some of the groups PSL$(2,q)$, where $q$ is a small prime power.
These results are partially new, and obtained through
both theoretical reasoning and application of computational techniques.  
}
\smallskip \hrule

\vspace*{.4in}

\begin{center}
{\large \bf 1 Introduction}
\end{center}

\bigskip

\noindent 
Let $\Sigma_g$ denote an orientable closed compact surface
of genus $g\geq 0$. Let $G$ be a finite group acting on $\Sigma_g$ preserving
orientation. It is known that there is a complex structure $(X, \Sigma_g)$
on the surface $\Sigma_g$ such that the group action can be realized 
as a subgroup of Aut$(X, \Sigma_g)$. Given such a surface $\Sigma_g$,
there possibly are uncountably many such complex structures $(X, \Sigma_g)$.
However, for a fixed $g \geq 2$, there are only finitely many groups 
of finite order which could possibly act faithfully on such surfaces.     
Indeed, it is well-known since the time of Hurwitz \cite{hur},
that for any $g \geq 2$ and finite subgroup $G$ of Aut$(X, \Sigma_g)$
we have $|G| \leq 84(g-1)$. 

\bigskip

\noindent 
It is however uninteresting to list the finite groups acting on surfaces
without the knowledge of the topological action. The following result, also 
known as Hurwitz's Theorem, essentially translates the topological language 
into combinatorial group theory. From now on we will consider only 
faithful actions, unless otherwise stated.

\bigskip

\noindent {\bf Theorem 1.1: (Riemann's Existence Theorem)}
A finite group $G$ acts on a surface $\Sigma_g$ of genus $g \geq 0$
if and only if 

\medskip

(1) there are integers $h \geq 0 $ and $n_1, \dotsc, n_r \geq 2$, for
some $r\geq 0$, such that the Riemann-Hurwitz relation 
\[
2(g-1) = |G|\left( 2(h-1) + \sum_{i=1}^{r} ( 1 - {\frac {1}{n_i}}) \right)
\]
holds, and 

\medskip

(2) there are elements $a_1, b_1, \dotsc, a_h, b_h, c_1, \dotsc ,c_r \in G$ 
such that
\[
G = \langle a_1, b_1, \dotsc, a_h, b_h, c_1, \dotsc c_r \rangle,
\]
fulfilling the long relation
\[ 
\prod_{i=1}^{h} [a_i, b_i]\cdot\prod_{j=1}^{r} c_j = 1  
\]
and such that
\[
|c_j| = n_j, \quad 1 \leq j \leq r  .
\]
\qed

\bigskip

\noindent 
Applying this to the genera $g = 0$ and $g = 1$
leads to the well-known result, indicated in Table \ref{tbl1},
on the finite groups which can act on such surfaces, see
for example \cite[Sect.I.3.4]{bre}. 

\begin{table}[here]
\begin{center}
\begin{tabular}{ | l | l | }
\hline
            genus & finite groups \\
\hline
\hline
        $0$ & $C_n$ (cyclic), $D_n$ (dihedral), $A_4$, $S_4$, $A_5$  \\
\hline
        $1$ & $(C_m \times C_n) : T$, where $T = C_2, C_3, C_4$ or $C_6$  \\
\hline
\end{tabular}
\end{center}
\caption{}\label{tbl1}
\end{table}


\noindent
Thus we will henceforth only consider actions of finite groups on
surfaces of genus $g\geq 2$. As was said earlier, for any fixed $g$ 
there are only finitely many groups in question. But changing
the point of view, and fixing a finite group, leads to the following notion:

\bigskip

\noindent {\bf Definition 1.2:} 
Let $G$ be a finite group. The set of all integers $g \geq 2$ which 
satisfies the conditions in Theorem $1.1$ is called the
{\bf genus spectrum} of $G$. We denote it by sp($G$). 

\bigskip

\noindent 
Let $\delta = (h; n_1, n_2, \dotsc, n_r)$ be a tuple of $r+1$ integers 
satisfying the conditions of Theorem $1.1$. 
We see that if 
$\delta^{\prime} = (h; n_1^{\prime}, n_2^{\prime}, \dotsc, n_r^{\prime})$
is another such tuple where 
$(n_1^{\prime}, n_2^{\prime}, \dotsc, n_r^{\prime})$ is a permutation of 
$(n_1, n_2, \dotsc, n_r)$, then $\delta^{\prime}$ also satisfies the
conditions of Theorem $1.1$.
Hence let ${\mathbb D}_r(G)$ be the quotient of the collection of all
tuples $\delta = (h; n_1, n_2, \dotsc, n_r)$ under the action of the
symmetric group $S_r$ on its last $r$ coordinates. 

\bigskip

\noindent 
We call the set 
${\mathbb D}(G) := \bigcup_{r \geq 0} {\mathbb D}_r(G)$ 
the {\bf data spectrum of $G$}, and its elements {\bf data} 
or {\bf signatures}.
Then there is the {\bf genus function}
{\bf g}$: {\mathbb D}(G) \longrightarrow {\mathbb N}_0$ defined by
\[ \delta = (h; n_1, n_2, \dotsc, n_r) \mapsto
{\mathrm \bf g}(\delta) := 
|G|\left( 
(h-1) + {\frac {1}{2}} \sum_{i=1}^{r} ( 1 - {\frac {1}{n_i}}) \right) + 1 
\]
sending a datum to its genus. Thus the image of {\bf g}, 
possibly after deleting a subset of $\{0,1\}$, 
is precisely the genus spectrum of $G$.
The basic arithmetic 
property of the genus spectrum sp($G$) is given as follows:

\bigskip

\noindent {\bf Theorem 1.2 {\cite{kul}}:} 
Let $G$ be a finite group. 
Then there exists a natural number $N_G$ such that: 

\medskip

(1) If $g \in$ sp($G$) then $g \equiv 1$ mod $N_G$, and

\medskip

(2) sp($G$) is cofinite in the set $1 + N_G {\mathbb N}$,
that is, up to finitely many exceptions for all $g$ satisfying
$g \equiv 1$ mod $N_G$ we have $g \in$ sp($G$).
\qed

\bigskip

\noindent 
The number $N_G$ for a finite group $G$ is called the
{\bf genus increment} for $G$, and is found 
easily as is outlined in Section $2$.

\bigskip

\noindent 
Let ${\mathcal FG}$ denote the class of all finite groups,
and ${\mathcal A}$ denote all sets of positive integers
which are cofinite subsets of some arithmetic set
of the form $1 + \lambda {\mathbb N}$, for some $\lambda \in {\mathbb N}$.
We consider the function 
{\bf {sp}} $: {\mathcal FG} \longrightarrow {\mathcal A}$,
which assigns to every finite group its genus spectrum, that is, 
{\bf {sp}}($G$) := sp($G$).
We immediately see that the function {\bf sp} is a decreasing function 
with respect to set-theoretic inclusion on both sides.
It is natural to ask the following questions:

\bigskip

\noindent {\bf Question 1:} Is {\bf sp} surjective?

\bigskip

\noindent 
It is difficult to answer this question since it is not at all easy 
to compute the spectrum of finite groups. However the set 
${\mathbb N}\setminus\{1\}$ is the spectrum of the cyclic groups $C_2$
and $C_3$, of order $2$ and $3$, respectively. Hence at least the full
set ${\mathbb N}\setminus\{1\}$ is in the image of {\bf sp}. 

\bigskip

\noindent 
The function {\bf sp} is clearly not injective, as seen in the 
previous paragraph.
However for a specific subclass $\mathcal C$ of ${\mathcal FG}$, we can
restrict the function {\bf sp} to $\mathcal C$, and ask if the restricted 
function {\bf sp}$\lvert_{\mathcal C}$ is injective. For example,
if ${\mathcal C}_p$ is the class of all finite cyclic $p$-groups,
with $p$ fixed, the restricted function {\bf sp}$\lvert_{{\mathcal C}_p}$ 
is known to be injective {\cite {kma}}. The following question is still open:

\bigskip

\noindent {\bf Question 2:}
Let $p$ be a fixed prime, and let ${\mathcal {AB}}_p$
be the class of all finite abelian $p$-groups. 
Is {\bf sp}$\lvert_{{\mathcal {AB}}_p}$ injective?

\bigskip

\noindent 
A partial answer is given in {\cite {tal}}, where it is shown that 
{\bf sp}$\lvert_{{\mathcal {AB}}_{(p,p^2)}}$ is injective, where 
${\mathcal {AB}}_{(p,p^2)}$ is the class of all
finite abelian $p$-groups of exponent up to $p^2$.

\bigskip

\noindent 
Coming back to the question of surjectivity of {\bf sp},
to get some idea we need a lot of information regarding the spectra 
of various types of finite groups. 
Towards this, we introduce more parameters:

\bigskip
  
\noindent {\bf Definition 1.3} 
Let $G$ be a finite group with genus increment $N_G$. 
The {\bf minimum genus} $\mu(G)$ and the {\bf stable upper genus} 
$\sigma(G)$ of $G$ are defined as follows:

\medskip

$\mu(G) := {\mathrm {min}} ~ {\mathrm {sp}}(G)$ 

\medskip

$\sigma(G) := {\mathrm {min}} 
\{ g \in {\mathrm {sp}}(G) ~:~  g + N_G{\mathbb N}_0 
\subseteq {\mathrm {sp}}(G) \}$  

\bigskip

\noindent 
Moreover, the finite set
\[
{\mathcal I}(G) := \{ g \geq \mu(G) : g \in 1 + N_G {\mathbb N}\} 
\setminus {\mathrm {sp}}(G) 
\]
is called the {\bf gap sequence} of $G$. 
Hence a complete spectrum sp($G$) of $G$ is determined by the quadruple
$(N_G, \mu(G), \sigma(G), {\mathcal I}(G))$. 

\bigskip

\noindent
There is a long series of articles written by various authors 
describing the minimum genus $\mu(G)$, also referred to as the 
{\bf (strong) symmetric genus}, for various finite groups $G$. 
These groups include cyclic groups {\cite {har}}, non-cyclic
abelian groups {\cite {mac}}, metacyclic groups {\cite {mic}},
alternating and symmetric groups {\cite {con1}}. 

\bigskip

\noindent
Moreover, the minimum genus is
known for all $26$ sporadic simple groups 
{\cite {con2, cww, lwi, wil1, wil2, wil3}}, as indicated
in Table \ref{tbl2}.
The corresponding actions are all given by some triangular datum 
of the form $(0; n_1, n_2, n_3) \equiv (n_1, n_2, n_3)$, where the
adjacent letters represent the numerals of the conjugacy classes 
the corresponding group generators belong to, following the notation 
of {\cite{ccnpw}}. For the groups Co$_1$ and Fi$_{24}^{\prime}$ 
we have not been able to figure out all pieces of information,
the uncertainties are indicated as question marks.

\begin{table}[t]
\begin{center}
\begin{tabular}{ | l | l | l |}
\hline
            group & minimum genus & action \\
\hline
\hline
        M$_{11}$ &  $631$ & $(2A, 4A, 11A)$ \\
\hline
        M$_{12}$ & $3,169$ & $(2B, 3B, 10A)$ \\
\hline
	J$_1$ & $2,091$ & $(2A, 3A, 7A)$ \\
\hline
	M$_{22}$ & $34,849$ & $(2A, 5A, 7A)$ \\
\hline
	J$_2$ & $7,201$ & $(2B, 3B, 7A)$ \\
\hline
	M$_{23}$ & $1,053,361$ & $(2A, 4A, 23A)$ \\
\hline
	HS & $1,680,001$ & $(2B, 3A, 11A)$ \\
\hline
	J$_3$ & $1,255,825$ & $(2A, 4A, 5A)$ \\
\hline
	M$_{24}$ & $10,200,961$ & $(3A, 3B, 4C)$ \\
\hline
	McL & $78,586,201$ & $(2A, 5A, 8A)$ \\ 
\hline
	He & $47,980,801$ & $(2B, 3B, 7D)$ \\
\hline
	Ru & $1,737, 216,001$ & $(2B, 3A, 7A)$ \\
\hline
	Suz & $11,208,637,441$ & $(2B, 4D, 5B)$ \\
\hline
	O'N & $9,600,323,041$ & $(2A, 3A, 8A)$ \\
\hline
	Co$_3$ & $5,901,984,001$ & $(2B, 3C, 7A)$ \\
\hline
	Co$_2$ & $1,602,478,080,001$ & $(2C, 3A, 11A)$ \\
\hline
	Fi$_{22}$ & $768,592,281,601$ & $(2C, 3D, 7A)$ \\
\hline
	HN & $3,250,368,000,001$ & $(2B, 3B, 7A)$ \\
\hline
	Ly & $616,252,131,000,001$ & $(2A, 3B, 7A)$ \\
\hline
	Th & $1,080,308,855,808,001$ & $(2A, 3C, 7A)$ \\
\hline
	Fi$_{23}$ & $85,197,301,526,937,601$ & $(2C, 3D, 8C)$ \\
\hline
	Co$_1$ & $86,620,350,136,320,001$ & $(2?, 3D, 8?)$ \\
\hline
	J$_4$ & $1,033,042,512,453,304,321$ & $(2B, 3A, 7A)$ \\
\hline
	Fi$_{24}^{\prime}$ & $14,942,925,109,412,639,539,201$ &$(2B,3E,7?)$ \\
\hline
	B & $86,557,947,525,550,545,649,532,928,$ & $(2D, 3B, 8x),$ \\
		& $000,001$ & $x \in \{ K, M, N \}$ \\
\hline
	M & $9,619,255,057,077,534,236,743,570,$ &  $(2B, 3B, 7B)$ \\
	    & $297,163,223,297,687,552,000,000,001$ &  \\
\hline
\end{tabular}
\end{center}
\caption{}\label{tbl2}
\end{table}

\bigskip

\noindent 
The minimum genus of a finite group corresponding to the datum
$(0; 2,3,7)$ is of a special importance: From the Hurwitz upper bound
we see that a possible genus $g \geq 2$ for a finite group $G$ 
would satisfy $g \geq {\frac {|G|}{84}} + 1$. Now, if $g$ corresponds 
to the datum $(0; 2, 3, 7)$ it follows from the 
Riemann-Hurwitz relation in Theorem $1.1$, that $g = {\frac {|G|}{84}} + 1$,
which is the lower bound for $g$. 
This gives rise to the question of determining finite groups $G$ which
have this property, known as {\bf Hurwitz groups}; for more details
the reader is suggested to look into the survey articles
{\cite {con3, con4}} and the references therein.   

\bigskip

\noindent 
The stable upper genus $\sigma(G)$ is rather harder to compute and only 
known for very few types of groups: for finite cyclic $p$-groups
{\cite {kma}}, metacyclic groups {\cite {wea}}, as well as $p$-groups of 
cyclic deficiency $\leq 2$ and elementary abelian $p$-groups for $p$ odd
{\cite {mta}}. The answer is also known partially for abelian $p$-groups 
{\cite {tal}}, and cyclic groups of order $pq$ for primes $p \neq q$
{\cite {owe}}. A description of genus spectra of finite $p$-groups of 
exponent $p$ is also given in {\cite {sar}}.



\bigskip

\noindent 
Encouraged by the similarities of the genus spectra for $p$-groups 
($p$ odd) of exponent $p$ whose lower central series have the same 
first two quotients, it was worth investigating the genus
spectra for the finite non-abelian $p$-groups whose lower central
series quotients would be relatively simple to study. The corresponding
groups are known as $p$-groups of {\bf maximal class} (or co-class $1$,
see {\cite{lmc}}). It is shown in {\cite {sar}} that if $p$ is an odd
prime and $n \leq p$ then all finite $p$-groups of maximal class and
of order $p^n$ would have the same genus spectrum.
The following question is however open:

\bigskip

\noindent {\bf Question 3:} 
Show that for $n$ large there are only two genus spectra for the 
finite $p$-groups of maximal class of order $p^n$ where $p$ is odd.

\bigskip

\noindent 
The story for $2$-groups of maximal class is relatively simpler. 
There are are only three isomorphism types of finite $2$-groups of
maximal class. The minimum and the stable upper genus for these 
groups, and hence the genus spectra, are all different. 
The proofs of this are outlined in Section $3$.

\bigskip 

\noindent 
There is rather a different and relatively complicated 
journey towards description of genus spectra for finite simple groups. 
The genus spectra of the first twelve sporadic simple groups 
(out of the total $26$ of them) have been calculated using
techniques of computational group theory and utilizing the
computer algebra system {\sf GAP} \cite{gap}. 
These result would encourage us to investigate the estimate
of the number of solutions of the associated diophantine equations 
which are not realizable in $G$, also known as bad solutions.
This is discussed in detail in Section $4$. 

\bigskip

\noindent {\bf Acknowledgement:} The first author would like to thank partial
funding credits to Lady Davis (2007-8) and Golda Meir (2008-9) 
post-doctoral scholarship funding (Hebrew University, Jerusalem),
and the hospitality of RWTH Aachen.

\vspace{.4in}


\begin{center}
{\large \bf 2 Arithmetic of the spectrum and groups of GK-type}
\end{center}

\bigskip

\noindent 
The genus increment $N_G$ of a finite group $G$ is 
connected to the $2$-groups of GK-type which we will briefly describe now.

\bigskip

\noindent {\bf Definition 2.1:} 
A finite $p$-group $G$ of exponent $p^e$ is said to be of
{\bf Gorenstein-Kulkarni type (GK-type)} if the set 
$\kappa(G) := \{ y \in G ~:~ |y| < p^e \}$ of elements 
of non-maximal order forms a subgroup of $G$ of index $p$.

\bigskip

\noindent {\bf Theorem 2.2 {\cite {kul}}:} 
Given a finite group $G$, for every prime $p$ let $G_p$ denote a 
Sylow $p$-subgroup of $G$. Denote by $p^{n_p}$ and $p^{e_p}$ 
the order and exponent of $G_p$, respectively. Then the genus increment
$N_G$ of $G$ is given as follows:
\[ 
N_G := \epsilon\cdot\prod_{p | |G|} p^{n_p - e_p}
\] 
where 
\[
\epsilon := \left\{ \begin{array}{ll}
              1 & \mbox{if $G_2=\{1\}$ or $G_2$ is of GK-type};\\
             1/2 & \mbox{if $G_2\neq\{1\}$ is not of GK-type}.
\end{array} \right.
\]
\qed

\bigskip

\noindent 
Although it is hard to find a finite group $G$, once a set
$S \in {\mathcal A}$ is given, which satisfies sp($G$) $= S$, 
it is relatively easy to find a group $G$ which satisfies 
$N_G = N$, whenever $N \in {\mathbb N}$ is given.

\bigskip

\noindent {\bf Proposition 2.3:} 
Let $N \in {\mathbb N}$. Then there exists a non-trivial finite 
group $G$, such that $N_G = N$.

\bigskip

\noindent {\bf Proof:} 
If $N = 1$, then let $G$ be any non-trivial cyclic group. 
Next suppose $N = p^{\alpha}$, where $\alpha \geq 1$ and $p$ is a prime.
If $p$ is odd, then for $G := C_{p^{\alpha}} \times C_{p^{\alpha}}$,
where $C_{p^{\alpha}}$ denotes the cyclic group of order $p^{\alpha}$,
we get $N_G = p^{\alpha}$; if $p=2$, then for 
$G := C_{p^{\alpha+1}} \times C_{p^{\alpha+1}}$
we get $N_G = 1/2 \cdot 2^{\alpha+1}=2^{\alpha}$.
Finally, if $N = p_1^{\alpha_1} p_2^{\alpha_2} \cdots p_k^{\alpha_k}$,
where the $p_i$ are pairwise distinct primes and 
$\alpha_i \geq 1$ for all $i$, let $G_i$ be the finite abelian group
with $N_{G_i} = p_i^{\alpha_i}$, as described in the last paragraph.
Then define $G$ to be the direct product $G := \prod_{i=1}^{k} G_i$
and note that $N_G = \prod_{i=1}^{k} N_{G_i}$.
\qed

\bigskip

\noindent 
For $p$ odd, the importance of finite $p$-groups of GK-type arrive
from the results of {\cite {mta}}, where groups having 
{\bf maximal exponent property (MEP)} are considered. 
It follows easily that groups with GK-type have MEP,
and hence the stable upper genus for
$2$-generated groups of GK-type is easily calculated, 
see {\cite[Thm.4.7]{mta}}.

\bigskip

\noindent {\bf Question 4:} 
Find a closed formula for the stable upper genus for the finite
$p$-groups of GK-type.

\bigskip

\noindent 
Looking at the structure of the finite $p$-groups
of GK-type, it would be interesting to classify them through certain
parameters. Towards this, we arrange them into trees as follows,
see {\cite {ms}}:

\bigskip

\noindent 
If $G$ is of GK-type, the subgroup $\kappa(G)$ of index $p$ might 
again turn out to be a $p$-group of GK-type. We inductively define 
$\kappa^{i+1}(G) := \kappa(\kappa^{i}(G))$, with $\kappa^{0}(G) = G$,
if $\kappa^{i}(G)$ is of GK-type. We call the subgroup $\kappa^{i}(G)$ 
the $i$-th {\bf GK-core} of $G$. We note that all subgroups 
$\kappa^{i}(G)$ of $G$, if defined, are characteristic and normal.
Conversely, $G$ would be referred to as a {\bf GK-extension} of 
$\kappa^{i}(G)$ of degree $|G/{\kappa^{i}(G)}|$. 

\bigskip

\noindent
This leads to the levelled directed graphs 
${\mathcal T}(R) = (V_{\mathcal T}, E_{\mathcal T})$ which are
defined as follows:
The set $V_{\mathcal T}$ of vertices consists of finite $p$-groups,
amongst which there is a unique group $R$,
also called the {\bf root} of ${\mathcal T}(R)$,  
which is not of GK-type, and is designated as the group of level $0$. 
Any other finite $p$-group 
$G$ occurring is a vertex of level $i>0$ if $\kappa^{i}(G) = R$. 
The set $E_{\mathcal T}$ of directed edges consists of the edges $e$ 
with initial vertex 
$H$ and terminal vertex 
$\kappa (H)$, for some $p$-group $H\in V_{\mathcal T}$ of GK-type. 

\bigskip

\noindent
Note that 
this implies that ${\mathcal T}(R)$ is a directed tree, that is,
a directed connected graph without circuits. Hence ${\mathcal T}(R)$
is called the {\bf GK-tree} associated with $R$.
The question arises how these trees look like. A partial answer
is given by the following statements:

\bigskip

\noindent {\bf Theorem 2.4 {\cite {ms}}:}
Let $R$ be a finite $p$-group of order $p^n$ and exponent $p^e$
which is not of GK-type. Then the tree ${\mathcal T}(R)$ is infinite 
if and only if exponent of the center $Z(R)$ of $R$  is $p^e$.
\qed

\bigskip

\noindent {\bf Definition 2.5:} 
A maximal infinite linear subgraph, that is, a subgraph topologically 
homeomorphic to the real line $\mathbb R$, of an infinite graph is 
called a {\bf stem}.

\bigskip

\noindent {\bf Theorem 2.6 \cite{ms}:} 
Let $R$ be a finite $p$-group of order $p^n$ and exponent $p^e$ which 
is not of GK-type, such that ${\mathcal T}(R)$ is infinite. Then a 
finite $p$-group $G$ of GK-type in ${\mathcal T}(R)$ lies on a stem
of ${\mathcal T}(R)$ if and only if there exists 
$t \in Z(G)\setminus \kappa^{1}(G)$ such that $G = \langle t, R \rangle$.


\medskip 

\noindent 
In this case, if $G$ lies on level $l>0$ of ${\mathcal T}(R)$, 
then $G$ has order $p^{n+l}$ and exponent $p^{e+l}$, and $|t| = p^{e+l}$.    
\qed

\bigskip

\noindent 
We proceed a little further along these lines:

\bigskip

\noindent {\bf Definition 2.7:}
A finite $p$-group $G$ is called {\bf regular}, if for every 
$x, y \in G$ we have $x^p y^p = (xy)^p c^p$ for some $c \in [H, H]$,
where $H = \langle x, y \rangle \leq G$ and
$[H, H]$ denotes the derived subgroup of $K$.

\bigskip

\noindent 
Hence in particular abelian $p$-groups, as well as $p$-groups
of exponent $p$ are regular.




\bigskip

\noindent {\bf Proposition 2.8:} 
Let $G$ be a finite $p$-group of GK-type which lies on a stem of an 
infinite GK-tree ${\mathcal T}(R)$. Then $G$ is regular if and only if 
all its GK-cores are regular, which holds if and only if $R$ is regular. 

\medskip

\noindent 
In particular, if $R$ is abelian or of exponent $p$,
then all GK-extensions of $R$
which lie on the stems of ${\mathcal T}(R)$ are regular. 

\bigskip

\noindent {\bf Proof:} Since any subgroup of a finite regular $p$-group 
is regular again, it is enough to show that if $R$ is regular and
$G$ is a GK-extension of $R$ such that exp($G$) $=$ exp($Z(G)$),
then $G$ is regular. 

\medskip

\noindent 
Now we have $t \in Z(G)$ such that $G = \langle t, R \rangle$,
and letting $x, y \in G$, we may write 
$x = t^{\alpha} a$ and $y = t^{\beta} b$ 
for some integers $\alpha, \beta$ and $a, b \in R$. Hence we have
$x^p y^p = t^{p(\alpha + \beta)} a^p b^p$
and $(xy)^p = t^{p(\alpha + \beta)} (ab)^p$.
Using the regularity of $R$, we have $a^p b^p = (ab)^p c^p$ 
for some $c \in [H,H]$ where $H = \langle a, b \rangle$. 
But $t$ centralizes $R$, which means $[H, H] = [K,K]$, 
where $K = \langle x, y \rangle$. Moreover, we get
\[
x^p y^p = t^{p(\alpha + \beta)} a^p b^p
= t^{p(\alpha + \beta)}(ab)^p c^p = (xy)^p c^p,
\]
which implies that $G$ is regular.
\qed

\bigskip

\noindent 
If we call these the {\bf regular central groups of GK-type},
denoted by ${\mathcal {RCGK}}_p$,
then we see that 
${\mathcal {AB}}_p \subseteq {\mathcal {RCGK}}_p$. Then we can ask:

\bigskip

\noindent {\bf Question 5:} Let $p$ be a fixed prime. 
Is {\bf sp}$\lvert_{{\mathcal {RCGK}}_p}$ injective?

\vspace{.4in}

\begin{center}
{\large \bf 3 Genus spectrum of $2$-groups of maximal class}
\end{center}

\bigskip

\noindent 
In this section we will show that the genus spectra of $2$-groups of maximal class are all different. Quite contrary to what one would be expecting in case of $p$-groups of maximal class where $p$ is odd.

\bigskip

\noindent For the terminologies specific to finite $p$-groups, we need to adopt the following specific notations. 

\bigskip

\noindent 
Let $p$ be a fixed odd prime and $G$ be a finite $p$-group of order $p^n$ and exponent $p^e$. Let $\delta = (h; n_1, \dotsc, n_r) \in {\mathbb D}(G)$. Since $n_i$ are powers of $p$ which could be repeated, we denote them by $\delta = (h; m_1, m_2, \dotsc, m_e)$ instead of $(h; p, \dotsc, p, p^2, \dotsc, p^2, \dotsc, p^e, \dotsc, p^e)$, where $p^i$ is repeated $m_i$ times in $\delta$. A set 
${\mathfrak G} = {\mathfrak G}_{\delta} 
= H_{\mathfrak G} \cup \bigcup_{i=1}^{e} E_{{\mathfrak G}, i})$ of generators corresponding to $\delta$ would be labelled as follows:

\smallskip

\noindent $\bullet$ Hyperbolic generators: 
$H_{\mathfrak G}= \{a_1, b_1, \dotsc, a_h, b_h\}$.

\smallskip 

\noindent $\bullet$ Elliptic generators of order $p^i$: 
$E_{{\mathfrak G}, i}= \{x_{i1}, x_{i2}, \dotsc, x_{i m_i}\} ~~ (1 \leq i \leq e)$, 

\smallskip

\noindent where the repeated elements are considered as different in the union $H_{\mathfrak G} \cup \bigcup_{i=1}^{e} E_{{\mathfrak G}, i}$.

\smallskip

\noindent The Riemann-Hurwitz equation takes the form:
\[
2(g-1) = p^n \left( 2(h-1) + \sum_{i=1}^{r} m_i ( 1 - {\frac {1}{p^i}} ) \right),
\] 
and the long relation looks like:
\[
\prod H_{\mathfrak G} \cdot\prod_{i=1}^{e} \prod E_{{\mathfrak G}, i} = 1,
\]
where 
\[
\prod H_{\mathfrak G} := \prod_{i=1}^{h} [a_i, b_i] ~~~ {\text{and}} ~~~
\prod E_{{\mathfrak G}, j} := \prod_{l=1}^{m_j} x_{j l}.
\]   

\bigskip

\noindent {\bf Data extensions 3.1:} This is a method which will be applied in the main proof. Let $\delta = (h; m_1, m_2, \dotsc, m_e) \in {\mathbb D}(G)$. Now every element $x$ of order $2^i$ in $G$ can be expressed as $x = x_1 x_2$, where $|x_1| = 2^j < 2^i$ and $|x_2| = 2^i$. Replacing $x$ by $x_1 x_2$ in the long relation and moving them to the right place using conjugates, we realize that $\delta^{\prime} = (h; m_1, m_2, \dotsc, m_{j-1}, m_j + 1, m_{j+1}, \dotsc, m_e) \in {\mathbb D}(G)$. Inductively, this shows that the multiplicities $m_i$ can be arbitrarily increased within ${\mathbb D}(G)$, if $m_i$ is not the last non-zero term in $\delta$ among the $m_k$'s. On the other hand arbitrarily many hyperbolic pairs of generators can be inserted on top of the existing ones. Hence if $m_t$ is the last non-zero term, $\delta^{\prime} = (h^{\prime}; m_1^{\prime}, \dotsc, m_{t-1}^{\prime}, m_t, 0, \dotsc, 0) \in {\mathbb D}(G)$ if $h^{\prime} \geq h, m_k^{\prime} \geq m_k$, where $1 \leq k \leq t-1$. 

\bigskip

\noindent To understand the stable upper genus of $2$-groups of exponent $2^e$, where $e \geq 5$, we need to look at the following diophantine equation:
\[
M = 2^e\cdot h + \sum_{i=1}^{e-1} m_i (2^{e-1} - 2^{e-1-i})
\]

\bigskip

\noindent Let $\Omega_e(2)$ denote the set of solutions $M$ of the above equation for which $h, m_i \geq 0$. For a fixed $e$, consider the $2$-adic expansion
\[
M = a_0 + a_1 2 + \dotsc + a_{e-1} 2^{e-1}
\]  
where $a_i = 0, 1$ for $0 \leq i \leq e-2$ and $a_{e-1} \geq 0$. Let $S_e (M)$ denote the sum $a_0 + a_1 + \dotsc + a_{e-1}$.

\bigskip

\noindent {\bf Theorem 3.2 {\cite {kma}}:}
We have
\[
\Omega_e(2) = \{ M \in {\mathbb N} : S_e (M) \geq e - 1 - i, 
a_i ~{\mathrm{is~the~first~non~zero~coefficient}}\}.
\]
\qed

\bigskip

\noindent Let $\delta_e (2)$ denote the smallest positive integer such that 
all $M \geq \delta_e (2)$ are realized as a solution of the above diophantine equation, for some $h, m_i \geq 0$. It is known from the properties of diophantine equations that such an integer $\delta_e (2)$ exists.

\bigskip

\noindent {\bf Theorem 3.3 {\cite {kma}}:} The least  stable solution $\delta_e (2)$ is given by:
\[
\delta_e (2) = (e-3) 2^{e-1} + 2
\]
\qed

\bigskip

\noindent There are three isomorphism types of $2$-groups of maximal class in each order $2^{e+1}$; namely, dihedral, quaternion and semi-dihedral $2$-groups,
being given in terms of presentations as follows, see {\cite{lmc}}:

\bigskip

\noindent (1) {\bf Dihedral Group} : $\langle y,x : x^2 = 1 = y^{2^e}, x^{-1}yx = y^{-1}
\rangle$

\bigskip

\noindent (2) {\bf Quaternion Group} : $\langle y,x : x^2 = y^{2^{e-1}}, x^4 = 1, x^{-1}yx =
y^{-1} \rangle$

\bigskip

\noindent (3) {\bf Semi-Dihedral Group} : $\langle y,x : x^2 = 1 = y^{2^e}, x^{-1}yx =
y^{2^{e-1}-1} \rangle$

\bigskip

\noindent {\bf Theorem 3.4:} 
Let $G$ be a finite $2$-group of order $2^{e+1}$, where $e \geq 5$, and of maximal class. Then the minimum and stable upper genus of $G$ are as given in 
Table \ref{tbl3}. 

\begin{table}[here]
\begin{center}
\begin{tabular}{ | l || l | l | }
\hline
            Type & Minimum Genus $\mu_0$ & Stable Upper Genus $\delta_0$ \\
\hline
\hline
        Dihedral & $2^{e-1}$ & $(e-3)2^e + 4$ \\
\hline
        Quaternion &  $2^{e-1}$ & $(e-3)2^e + 2^{e-1} + 3$ \\
\hline
        Semi-Dihedral & $2^{e-2}$ & $(e-3)2^e + 2^{e-2} + 3$ \\
\hline

\end{tabular}
\end{center}
\caption{}\label{tbl3}
\end{table}


\noindent All these three groups of order $2^{e+1}$ contain a (maximal) cyclic subgroup of index $2$. Since none of these are of GK-type, the genus increment $N_G$ for all of them equals $1$. We will only show the details of the remaining calculations for dihedral groups. The calculations for other two groups are similar. The proof for the dihedral group relies on the following lemmas:

\bigskip

\noindent {\bf Lemma 3.5:} Let $G$ be a dihedral group of order $2^{e+1}$ with generators $x, y$ satisfying the relations given above. If $G = \langle xy^a, xy^b \rangle$, then $a-b$ is odd. If $G = \langle y^a, xy^b \rangle$, then $a$ is odd.

\bigskip

\noindent {\bf Proof of Lemma 3.5:} The Frattini subgroup is $\Phi(G) = \langle y^2 \rangle$, and hence the quotient equals $G/{\Phi(G)} = C_2 \times C_2 = \langle {\overline{x}}, {\overline{y}} \rangle$. Now if $a-b$ is even, we have ${\overline{xy^a}} = {\overline{xy^b}}$, which implies $G/{\Phi(G)}$ is cyclic, a contradiction. In the second case if $a$ is even, it again implies $G/{\Phi(G)}$ is cyclic, a contradiction. \qed

\bigskip

\noindent The elements of order $2$ in the dihedral group $G$ which are of the form $xy^a$, for $a \in {\mathbb Z}$, are often referred to as reflections in geometry. Next we note that:

\bigskip

\noindent {\bf Lemma 3.6:} The number of reflections among the elliptic generators $E_{{\mathfrak G},1}$ (with repetitions allowed) is even.

\bigskip

\noindent {\bf Proof of Lemma 3.6:} Assume that $E_{{\mathfrak G},1}$ contains $k$ reflections, say $xy^{i_1}, xy^{i_2}, \dotsc, xy^{i_k}$. Since the reflections in $G$ form a union of conjugacy classes, we note that $\prod E_{{\mathfrak G},1} = xy^{i_1} xy^{i_2} \dotsc xy^{i_k} \xi$ for some $\xi \in \{ 1, y^{2^{e-1}} \}$. Since $xy^a xy^b = y^{b-a}$, for $k$ odd we get $\prod E_{{\mathfrak G},1} = xy^j$ for some $j \in {\mathbb Z}$. Now since $\prod H_{\mathfrak G}, \prod E_{{\mathfrak G}, i} \in \langle y \rangle$, where $i \geq 2$, using the long relation we get $xy^j \in \langle y \rangle$, a contradiction.
\qed

\bigskip

\noindent {\bf Lemma 3.7:} If $\delta = (h; m_1, \dotsc, m_e) \in {\mathbb D}(G)$ with $h \leq 1$ and $m_{e-1} = m_e = 0$, then $E_{{\mathfrak G},1}$ contains at least two reflections.

\bigskip

\noindent {\bf Proof of Lemma 3.7:} Assume the contrary, then we have $E_{{\mathfrak G},i} \subseteq \Phi(G) = \langle y^4 \rangle$ for all $i$. Hence $H_{{\mathfrak G},i}$ generates $G$, thus $h =1$ and this contradicts the long relation using Lemma $3.6$.
\qed  

\bigskip

\noindent {\bf Proof of Theorem 3.1} for the dihedral group $G$ of order $2^{e+1}$: We will first show that the minimum genus $\mu(G)$ corresponds to the datum $\delta_0 = (0; 3, 0, \dotsc, 0, 1)$ given by {\bf g}$(\delta_0) = 2^{e-1}$. 

\bigskip

\noindent Let $\delta = (h; m_1, \dotsc, m_e) \in {\mathbb D}(G)$ correspond 
$\mu$. If $h \geq 2$, then {\bf g}$(\delta) \geq 2^{e+1}+1$. Hence $h = 0$ or $1$.

\bigskip

\noindent We consider the case $h = 1$. If $m_e \neq 0$, 
since $\gamma_2(G) = \langle y^2 \rangle$ we have from the long relation that
$m_e \geq 2$. Then {\bf g}$(\delta) \geq 2^{e+1} (1 - {\frac {1}{2^e}}) + 1 > 2^{e-1}$, hence $m_e = 0$. Next, if $m_{e-1} \geq 1$, then {\bf g}$(\delta) \geq 2^e (1 - {\frac {1}{2^{e-1}}}) + 1 > 2^{e-1}$, hence $m_{e-1} = 0$.

\bigskip

\noindent By lemma $3.4$ we have that $E_{{\mathfrak G},1}$ contains at least two reflections, that is, $m_2 \geq 2$. Hence {\bf g}$(\delta) \geq 2^e + 1 > 2^{e-1}$. This proves that $h \neq 1$.

\bigskip 

\noindent Now we consider the case $h=0$. Since the group cannot be generated without at least one reflection, we have $m_2 \geq 2$. If $m_2 = 2$ and $m_e \geq 2$, then g$(\delta) \geq 2^e$. But with $m_2 = 2$ and $m_e = 1$, we need some more elliptic generators to make the right side of the Riemann-Hurwitz formula positive, say $m_i \geq 1$. Then g$(\delta) \geq 2^{e+1} ({\frac {3}{2}} - {\frac {1}{2^i}} - {\frac {1}{2^e}}) + 1 > 2^e$. Hence with $m_2 = 2$ we have $m_e = 0$. But the product of two reflections generating $G$ is $2^e$. Hence if $m_2 \geq 2$ and $m_e = 0$ we would need at least four reflections in $E_{{\mathfrak G},1}$. Hence $m_2 \geq 4$. Now with $m_2 \geq 5$ we have g$(\delta) \geq 2^{e-1} + 1$. Hence we consider $m_2 = 4$. We would then need more elliptic generators to make the genus $> 1$, say $m_i \geq 1$ for some $i$. Then g$(\delta) \geq 1 + 2^e - 2^{e-i} \geq 2^{e-1} + 1$.  

\bigskip

\noindent This shows that the minimum genus of $G$ is $\mu(G)=2^{e-1}$. We will next show that the stable upper genus $\sigma(G)$ equals $(e-3)2^e + 4$. For this we show, recalling hat $N_G = 1$:

\smallskip

\noindent (1) If $g \geq (e-3)2^e + 4$, then $g \in$ sp$(G)$. \\
\noindent (2) $(e-3)2^e + 3$ is not realized by any datum of $G$.

\bigskip

\noindent {\bf Proof of (1):} For such an integer $g$, consider $V_g := g + 2^{e+1} - 1$. Equivalently we will show that if $V_g \geq (e-1) 2^e + 3$, then there exists a datum $\delta = (h; m_1, m_2, \dotsc, m_e) \in {\mathbb D}(G)$ such that {\bf g}$(\delta) = g$. 

\bigskip

\noindent Case I : $V_g$ is odd.

\smallskip

\noindent This ensures that $m_e$ is odd. We first consider $m_e = 1$, and then any $\delta^{\prime} = (h; m_1, m_2, \dotsc, m_{e-1}, 1)$ is a datum if $m_1 \geq 2$, since it is a proper extension of $(0; 2, 0, \dotsc, 0, 1)$, where $(0; 2, 0, \dotsc, 0, 1)$ is realized by ${\mathfrak G} = \{ x, xy, y^{-1} \}$ and does not associate to a genus $\geq 2$. Now
\[
{\frac {V_g - 2^e + 1}{2}} = 2^e \cdot h + \sum_{i=1}^{e-1} (2^{e-1} - 2^{e-1-i}) m_i.
\]
From Theorem $3.3$ it follows that all odd integers $V_g$ with ${\frac {V_g - 2^e + 1}{2}} \geq (e-3) 2^{e-1} + 2 + 2 . (2^{e-1} - 2^{e-2})$ can be realized. This means that all odd integers $V_g \geq (e-1) 2^e + 3$ with $m_e = 1$ can be realized. For an arbitrary odd $m_e$, insert a slot $u u^{-1}$ in the elliptic part of the long relation, where the order of $u$ is $2^e$. Finally, this means that all even integers $g \geq (e-3) 2^e + 4$ are genera.

\bigskip

\noindent Case II : $V_g$ is even.

\smallskip

\noindent Here $m_e$ is even. Consider $m_e = 0$, and then any tuple $\delta^{\prime} =$ \linebreak $(h; m_1, m_2, \dotsc, m_{e-1}, 0)$ is a datum if $m_1 \geq 4$, since it is a proper extension of $(0; 4, 0, \dotsc, 0)$, where $(0; 4, 0, \dotsc, 0)$ is realized by ${\mathfrak G} = \{ x, xy, x, xy^{-1} \}$ and does not associate to a genus $\geq 2$. Now
\[
{\frac {V_g}{2}} = 2^e \cdot h + \sum_{i=1}^{e-1} (2^{e-1} - 2^{e-1-i}) m_i.
\]
From Theorem $3.3$ it follows that all odd integers $V_g$ with ${\frac {V_g}{2}} \geq (e-3) 2^{e-1} + 2 + 4. (2^{e-1} - 2^{e-2})$ associate a realized data. This mean all even integers $V_g \geq (e-1) 2^e + 4$ realize a data, and therefore all odd integers $g \geq (e-3) 2^e + 5$ are genus.
\qed

\bigskip

\noindent {\bf Proof of (2):} In this part we will examine all possible solutions of the equation:
\[
V_g =  (e-1)\cdot 2^e + 2 = 2^{e+1} \cdot h + \sum_{i=1}^{e} (2^e - 2^{e-i}) m_i
\]

\bigskip

\noindent {\bf Lemma 3.8:} $m_e$ is even and $\not\geq 4$.

\smallskip

\noindent {\bf Proof of Lemma 3.8:} It is clear that $m_e$ is even. Now from the above equation we get 
\begin{equation*}\tag{A}\label{eqa}
2^e \cdot (2h + m_1 + \dotsc + m_e) - M(m_1, \dotsc, m_e) = (e-1) 2^e + 2 
\end{equation*}
where $M(m_1, \dotsc, m_e) = 2^{e-1} m_1 + 2^{e-2} m_2 + \dotsc + 2 m_{e-1} + m_e$. 

\smallskip

\noindent This implies that $M(m_1, \dotsc, m_e) \equiv -2$ mod $2^e$.

\smallskip

\noindent Sub-case I : $m_e \equiv 2$ mod $2^2$. It is possible to define the following equations, with integers $r_i \geq 0$:
\begin{eqnarray*}
m_e =& 2 + 2^2 r_e \\
m_{e-1} + 2 + 2 r_e =& 2 r_{e-1} \\
m_{e-2} + r_{e-1} =& 2 r_{e-2} \\
m_{e-3} + r_{e-2} =& 2 r_{e-3} \\
m_{e-4} + r_{e-3} =& 2 r_{e-4} \\
\dotsc & \\
m_2 + r_3 =& 2 r_2 \\
m_1 + r_2 =& 2 r_1
\end{eqnarray*}


\noindent From equation (\ref{eqa}) we get 
\[
2^e \cdot 2h = 2^e (e-1 + r_1 - (m_1 + \dotsc + m_e)).
\]
Since $h \geq 0$ we have that $m_1 + \dotsc + m_e - r_1 \leq e-1$, which implies
\begin{equation*}\tag{B}\label{eqb}
r_1 + r_2 + \dotsc + r_{e-1} + 2 r_e \leq e-1.
\end{equation*}
Now suppose $m_e \geq 6$, which implies $r_e \geq 1$. From the list of relations we get that $m_{e-1} \geq 2$, and $r_i \geq 1$ if $1 \leq i \leq e-2$. This contradicts (\ref{eqb}). Hence $r_e = 0$ and consequently $m_e = 2$. 

\smallskip

\noindent Sub-case II : $m_e \equiv 0$ mod $2^2$. It is possible to define the
following equations, with integers $r_i^{\prime} \geq 0$:
\begin{eqnarray*}
m_e =& 2^2 r_e^{\prime} \\
m_{e-1} + 1 + 2 r_e^{\prime} =& 2 r_{e-1}^{\prime} \\
m_{e-2} + r_{e-1}^{\prime} =& 2 r_{e-2}^{\prime} \\
m_{e-3} + r_{e-2}^{\prime} =& 2 r_{e-3}^{\prime} \\
m_{e-4} + r_{e-3}^{\prime} =& 2 r_{e-4}^{\prime} \\
\dotsc& \\
m_2 + r_3^{\prime} =& 2 r_2^{\prime} \\
m_1 + r_2^{\prime} =& 2 r_1^{\prime}
\end{eqnarray*}



\noindent From equation (\ref{eqa}) we get
\[
2^e \cdot 2h = 2^e (e - (r_1^{\prime} + r_2^{\prime} + \dotsc + r_{e-1}^{\prime} + 2 r_e^{\prime})).
\]
Since $h \geq 0$ we get
\begin{equation*}\tag{C}\label{eqc}
r_1^{\prime} + r_2^{\prime} + r_{e-1}^{\prime} + 2 r_e^{\prime} \leq e.
\end{equation*}
Now if $m_e \geq 4$ we have $r_e^{\prime} \geq 1$. Also inductively we obtain $r_{e-2}^{\prime} \geq 2$, and $r_i^{\prime} \geq 1$ if $1 \leq i \leq e-3$. This contradicts (\ref{eqc}). Hence $m_e = 0$.
\qed

\bigskip

\noindent {\bf Lemma 3.9:} If $m_e = 2$, then the only solution of $(A)$ is given by 
\[
h = 0, m_i = 1 ~~\text{for}~~ 1 \leq i \leq e-2, m_{e-1} = 0.
\]

\smallskip

\noindent {\bf Proof of Lemma 3.9:} From the list of equations in Lemma $3.8$ (Sub-case I) it follows that, if $m_{e-1} \geq 1$, then $r_{e-2} \geq 2$ and for the other ones $r_i \geq 1$, which contradicts equation (\ref{eqb}). If $m_{e-1} = 0$, then $r_i \geq 1$ for $1 \leq i \leq e-1$. Hence from (\ref{eqb}), the only choice is $r_i = 1$ for $1 \leq i \leq e-1$. This gives the result from the list of equations given and $\ref{eqa}$.
\qed

\bigskip

\noindent {\bf Lemma 3.10:} If $m_e = 0$, then $m_1 \leq 3$ and $h = 0$.

\bigskip

\noindent {\bf Proof of Lemma 3.10:} If $m_e = 0$, then the similar analysis yields $r_i^{\prime} \geq 1$ for $1 \leq i \leq e-1$, and we get the relation
\[
e-1 \leq r_1^{\prime} + r_2^{\prime} + \dotsc + r_{e-1}^{\prime}. 
\]
If all of $r_i^{\prime} = 1$, then from equation $(A)$ we get that $h = 1/2$, a contradiction. Hence $r_i^{\prime} = 1$ for all $i$, except one factor, say $r_j = 2$ for some $1 \leq j \leq e-1$. Then we obtain:
\[
m_e = 0, m_{e-1} = 1, \dotsc, m_{j+1} = 1,
\]
\[
m_j = 3, m_{j-1} = 0, m_{j-2} = 1, \dotsc, m_1 = 1,
\]
with the exceptional solutions
\[
m_e = 0, m_{e-1} = 3, m_{e-2} = 0, m_{e-3} = 1, \dotsc, m_1 = 1
\]
and
\[
m_e = 0, m_{e-1} = 1, \dotsc, m_2 = 1, m_1 = 3.
\]
In all these cases, we have $h = 0$.
\qed

\bigskip

\noindent {\bf Back to proof of (2):} If $m_e = 2$, then $m_1 = 1, h = 0$ and hence any set $\mathfrak G$ with this datum can generate at most the subgroup $\langle y \rangle$. If $m_e = 0$, we have $h = 0$, and hence we need to generate $G$ by the reflections in $E_{{\mathfrak G}, 1}$, which cannot be more than two, hence by Lemma $3.6$ are exactly two. Then $\prod E_{{\mathfrak G}, 1}$ is an odd power of $y$, whereas $\prod_{i \neq 1} \prod E_{{\mathfrak G}, i}$ is an even power of $y$. This contradicts the long relation. 
\qed

\vspace{.4in}

\begin{center}
{\large \bf 4 Computation of genus spectra of finite simple groups}
\end{center}

\bigskip

\noindent 
We describe how the genus spectrum of a given finite group $G$
can be determined explicitly using tools from computational group
theory. To this end, let $\Delta(G)$ denote the set of non-trivial 
{\bf periods} of $G$, that is,
$\Delta(G) := \{ |x| \in {\mathbb N} : 1\neq x \in G \}
 =\{n_1,\dotsc, n_k\}$, say, for some $k\in{\mathbb N}_0$.
Then the {\bf Auxiliary Euler-Characteristic Equation (AEC equation)}
is defined by 
\[
2(g-1) = |G| \left( 
2(h-1) + \sum_{n \in \Delta(G)} a_n (1 - {\frac {1}{n}}) \right),
\]
where $g, h \in {\mathbb N}_0$, and $a_n \in {\mathbb N}_0$ 
for all $n \in \Delta(G)$. 
Letting $N_G$ still being the genus increment associated with $G$, 
the associated reduced AEC equation is:
\[
{\tilde {g}} := {\frac {(g-1)}{N_G}} = 
{\frac {|G|}{N_G}} \left( 
(h-1) + {\frac {1}{2}} \sum_{n \in \Delta(G)} a_n (1 - {\frac {1}{n}}) \right),
\]
where $\tilde {g}$ is called the {\bf reduced genus} 
corresponding to the genus $g$.

\bigskip

\noindent 
Let $D(G) := \bigcup_{\tilde {g}\geq 1}D_{\tilde {g}}(G)$, where for 
a reduced genus $\tilde {g}\geq 1$ the set $D_{\tilde {g}}(G)$
denotes the set of all solutions $\delta$ of the reduced AEC equation,
written as $\delta := (h; n_1^{[a_1]}, \dotsc, n_k^{[a_k]})$.
If $\delta$ is realized in the group $G$ as described by 
Riemann's Existence Theorem we notice that the tuple $\delta$  
just is a datum of $G$. This leads to the following:

\bigskip

\noindent {\bf Definition 4.1:} 
A solution $\delta\in D_{\tilde {g}}(G)$ 
to the reduced AEC equation of $G$ is called {\bf bad} if
$\tilde {g} \geq \tilde {\mu} (G):=\frac {\mu(G)-1}{N_G}$,
where the latter number is called the {\bf minimum reduced genus} of $G$,
and $\delta \not\in {\mathbb D}(G)$, that is $\delta$ is not a datum.

\bigskip

\noindent 
Thus letting $\tilde {\sigma} (G):=\frac {\sigma(G)-1}{N_G}$ be the
{\bf stable upper reduced genus} of $G$, describing the
spectrum of $G$ completely is equivalent to determining
$\tilde {\mu} (G)$ and $\tilde {\sigma} (G)$ as well as 
$D(G)$ and the bad solutions. Note that the bad solutions 
necessarily belong to the interval 
$\tilde {\mu} (G) \leq \tilde {g} \leq \tilde {\sigma} (G)$.  

\bigskip

\noindent 
Now, to determine $\tilde {\sigma} (G)$ we need to look at 
a diophantine problem due to Frobenius popularly known as the 
"Coin Problem": 
Given positive integers $n_1, n_2, \dotsc, n_k$ having no non-trivial 
common divisor, consider the linear equation
\[
n = n_1 x_1 + n_2 x_2 + \dotsc + n_k x_k,
\]
and call a non-negative integer $n$ dependent on $n_1, n_2, \dotsc, n_k$ if
there exists a solution $(x_1, x_2, \dotsc, x_k)$ of
the above equation in the non-negative integers,
and independent from $n_1, n_2, \dotsc, n_k$ otherwise.
It is known that there is a bound $f$ such that all $n$ exceeding
$f$ are dependent, and the largest independent integer
$f=f(n_1, n_2, \dotsc, n_k)$ 
is called the Frobenius number of $n_1, n_2, \dotsc, n_k$. 

\bigskip

\noindent 
There are no formulae known for $f(n_1, n_2, \dotsc, n_k)$ whenever 
$k \geq 4$; see {\cite {bsh, rod}}. But at least it is known that
$f(n_1, n_2) = (n_1 - 1)(n_2 - 1) - 1$. 
Our main interest is in the recurrence formula in case there are 
common divisors for any $(k-1)$ of the given numbers, which occurs 
always for the AEC equations.

\bigskip

\noindent 
{\bf Theorem 4.2 {\cite{bsh}}:} 
If $d$ divides $n_2, n_3, \dotsc, n_k$, then
\[
f(n_1, n_2, \dotsc, n_k) 
= d \cdot f(n_1, {\frac {n_2}{d}}, \dotsc, {\frac {n_k}{d}}) + n_1 (d-1).
\]
\qed

\bigskip

\noindent 
The existence of a datum $(0; n_1, n_2, n_3)\in {\mathbb D}_3(G)$ 
can be explicitly computed by the following character sum formula,
where ${\mathrm {Irr}}(G)$ denotes the set of irreducible complex
characters of $G$:

\bigskip

\noindent 
{\bf Theorem 4.3}, see {\cite[Thm.2.12]{gor}}:
Let $G$ be a finite group, and let $C_1, C_2, \dotsc, C_k$ 
denote the conjugacy classes in $G$, with representatives 
$g_1, g_2, \dotsc, g_k$.
Then the number of solutions to the equation $xyz = 1$ in $G$ 
with $x \in C_r$ and $y \in C_s$ and $z\in C_t$ is given by
\[
c_{C_r, C_s, C_t} = {\frac {|C_r|.|C_s|.|C_t|}{|G|}}\cdot 
\sum_{\chi \in {\mathrm {Irr}}(G)}
{\frac {\chi(g_r) \chi(g_s) \chi(g_t)}{\chi(1)}}.
\]
\qed

\bigskip

\noindent 
The number $c_{C_r, C_s, C_t}$ 
is also referred to as the associated 
{\bf class multiplication coefficient}.

\bigskip

\noindent 
The data extensions are slightly different from what we have done in the
previous section. Let $(h; n_1, n_2, \dotsc, n_r)\in {\mathbb D}_r(G)$ 
denote a datum of $G$ 
corresponding to the genus $g \geq 2$, allowing repetitions of the 
periods $n_i$, with a possible generating sequence 
$\{ a_1, b_1, \dotsc, a_h, b_h; c_1, \dotsc, c_r \}$
as in Riemann's Existence Theorem.
Let $c_j = c'_1 c'_2$ where the $c'_i$ are non-trivial elements of $G$ 
of order $n'_i$. Then replacing $c_j$ in the long relation by $c'_1c'_2$ 
we see that $(h; n_1, n_2, \dotsc, n_{j-1}, n'_1, n'_2, n_{j+1}, \dotsc, n_r)$ 
is a datum in ${\mathbb D}_{r+1}(G)$ corresponding to the genus
$g^{\prime}$, where 
it is straightforward to check that $g^{\prime} \geq g$.

\bigskip

\noindent
We now turn to explicit calculations for some finite simple groups. 
First we consider the infinite family PSL$(2,q)$, where $q$ is a 
prime power; recall that these groups are simple if and only if $q\geq 4$.

\bigskip

\noindent {\bf Theorem 4.4}, see \cite[Satz.II.8.10]{hup}: 
Let $G =$ PSL$(2, q)$, where $q=p^f$ is a prime power, and let
$l$ be a prime dividing $|G|$. Then we have:

\medskip

(1) If $2\neq l\neq p$, then the Sylow $l$-subgroups of $G$ are cyclic. 
 
\medskip

(2) If $p\neq 2$, then the Sylow $2$-subgroups of $G$ are dihedral.

\medskip

(3) The Sylow $p$-subgroups are elementary
abelian of order $q$.
\qed

\bigskip

\noindent
Thus, the genus increment of PSL$(2, q)$,
for any odd prime power $q=p^f$, equals $q/p=q^{f-1}$,
while the genus increment of PSL$(2, 2^f)$ equals
$2^{f-2}$ for $f\geq 2$, and that of PSL$(2, 2)$ equals $1$.
In particular, the genus increment of PSL$(2, q)$,
where $p$ is any prime, equals $1$.

\bigskip

\noindent
As for the minimum genus of PSL$(2, p)$, where $p$ is a prime,
we have the following results.
Recall that PSL$(2, 4)\cong$ PSL$(2, 5)\cong A_5$ acts
on the Riemann sphere, so its genuine minimum genus indeed 
equals $0$, as indicated below, while in our considerations 
of genus spectra we only look at genera $g\geq 2$.

\bigskip

\noindent {\bf Theorem 4.5 {\cite{gsj}}:}
For any prime $p \geq 13$ let
\[
d_p := {\mathrm {min}} \left\{ e : e \geq 7, 
~{\mathrm {either}}~ e | {\frac {p-1}{2}}
~{\mathrm {or}}~ e | {\frac {p+1}{2}} \right\}.
\]
Then the minimum genus of PSL$(2, p)$, where $p\geq 5$ is a prime,
comes from the following data:

\medskip

(1) $(0;2, 3, p)$ if $p = 5, 7, 11$,

\medskip

(2) $(0;2, 5, 5)$ 
if $p \geq 13,\, p \equiv \pm 1$ mod $5,\, p \not\equiv \pm 1$ 
mod $8$ and $d_p \geq 15$, 

\medskip

(3) $(0;3, 3, 4)$
if $p \geq 13,\, p \not\equiv \pm 1$ mod $5,\, p \equiv \pm 1$
mod $8$ and $d_p \geq 12$, 

\medskip

(4) $(0;2, 4, 5)$ 
if $p \geq 13,\, p \equiv \pm 1$ mod $5,\, p \equiv \pm 1$ 
mod $8$ and $d_p \geq 9$, 

\medskip

(5) $(0;2, 3, d_p)$ in all other cases.
\qed

\bigskip

\noindent
More generally, the minimum genus for PSL$(2, q)$, where $q$ is
an arbitrary prime power, has been determined in \cite{gsj2}. 
Moreover, in \cite{voo} there covering groups SL$(2,q)$ have been treated.
Thus we ask for the minimum stable genus and the bad solutions.
If we restrict ourselves again to the case $p$ a prime, we
have the following result:

\bigskip

\noindent {\bf Theorem 4.6 {\cite[Cor.5.10]{ost}}:}
Let $G=$ PSL$(2, p)$, where $p\geq 13$ is a prime.
Then for the data spectrum of $G$ we have 
${\mathbb D}(G)=D(G)$, in particular there are no bad solutions.
\qed

\bigskip

\noindent 
Thus, the spectrum of $G$ can be just read off from the 
solutions of the AEC equation:
Letting $\Delta(G)=\{n_1,\dotsc,n_k\}$ be the set of periods of
$G$, then the stable upper genus of $G$ is given as
$$ \sigma(G)=2+|G|+f(|G|(n_1-1)/2n_1,\dotsc,|G|(n_k-1)/2n_k) ,$$
where $f$ denotes the associated Frobenius number.

\bigskip

\noindent
Actually, it can be checked explicitly that the assertion
of Theorem 4.5 also holds for $p=5,7,11$. Indeed, we take 
the opportunity to present some details for the case $p=7$,
in order to indicate the flavor of the computations, and
to show how techniques from computational group theory, 
using the computer algebra system {\sf GAP} \cite{gap} and
its character table library \cite{CTblLib}, come into play:


\bigskip

\noindent 
The group $G=$ PSL$(2, 7)$ has order $168$ and periods
$\Delta(G) = \{ 2, 3, 4, 7 \}$. From Theorem $4.4$ we recover
the well-known fact that the minimum genus of $G$ equals $3$,
and corresponds to the datum $(0; 2, 3, 7)$, establishing $G$ 
as a Hurwitz group. The reduced AEC of $G$ is
\[
g-1={\tilde {g}} = 168(h-1) + 42 a_2 + 56 a_3 + 63 a_4 + 72 a_7,
\]
The associated Frobenius number being $f(168, 42, 56, 63, 72) = 565$,
we conclude that ${\tilde {g}} = 565-168=397$ is not a reduced genus.

\bigskip

\noindent 
Since {\sf GAP} \cite{gap} shows that for the class
multiplication coefficients in $G$ we $c_{7, 7, i} \neq 0$,
for all $i \in \Delta(G)$, the extension principle shows that
all tuples $(h; 2^{[a_2]}, 3^{[a_3]}, 4^{[a_4]}, 7^{[a_7]})$ 
with $h \geq 0, a_2 \geq 1, a_3 \geq 1, a_4 \geq 0, a_7 \geq 1$
are actually data of $G$.
Thus all ${\tilde {g}}>565-168+(42 + 56 + 72)=567$
are reduced genera of $G$. We will prove that all integers ${\tilde {g}}$
with $398 \leq {\tilde {g}} \leq 567$ are genera of $G$, and thus in
particular are solutions of the reduced AEC, implying that 
$\sigma(G)-1={\tilde {\sigma}}(G)=398$. Moreover, we will prove
that all solutions of the reduced AEC come from data of $G$,
so that there are no bad solution. 

\bigskip

\noindent 
Using {\sf GAP} \cite{gap} we verify that the above assertions follow
from the extension principle, as soon as we have shown that the
following actually are data: (Indeed, we cannot do with less.)
\[
(0; 2, 4, 7), 
(0; 2, 7, 7), 
(0; 3, 3, 4), 
(0; 3, 3, 7), 
(0; 3, 4, 4), 
\]
\[
(0; 3, 4, 7), 
(0; 3, 7, 7), 
(0; 4, 4, 4), 
(0; 4, 7, 7), 
(0; 7, 7, 7).
\]

\bigskip

\noindent 
To this end, firstly {\sf GAP} \cite{gap} shows that all the
relevant class multiplication coefficients are non-zero. Hence
for any case there are triples of elements fulfilling the  
associated long relation. Next, recall that, by {\cite {ccnpw}},
$G$ has precisely our conjugacy classes of maximal subgroups,
these consist of groups of shape $S_4$ (order $24$), $7 : 3$ (order $21$), 
$D_8$ (order $8$), $D_6$ (order $6$), respectively.
In particular, the only maximal subgroup having order divisible
by $7$ is of shape $7 : 3$, which moreover has a normal Sylow $7$-subgroup.
Hence triples corresponding to $(0; 2, 4, 7)$, $(0; 2, 7, 7)$,
$(0; 3, 4, 7)$, $(0; 3, 7, 7)$ and $(0; 4, 7, 7)$ necessarily generate $G$.

\bigskip

\noindent 
Moreover, from $c_{7X,7X,7X} = c_{7Y,7X,7X} = c_{7X,7Y,7X}=24$
and $c_{7Y,7Y,7X} = 216$, where $\{X,Y\}=\{A,B\}$, we conclude
that there are $576$ triples associated with $(0; 7, 7, 7)$.
If such a triple does not generate $G$, then it necessarily
generates a cyclic group of order $7$. But there are precisely
eight such subgroups in PSL$(2, 7)$, each of which contains
$30$ such triples. Hence there are generating triples of this shape.
Similarly, there are $1008$ triples in $G$ associated with $(0;3,3,7)$,
but in any of the eight subgroups isomorphic to $7 : 3$ there are
precisely $84$ such triples; hence there are triples of this shape
generating $G$. 

\bigskip

\noindent
Next we consider $(0; 3, 3, 4)$: Assume a corresponding triple does not
generate $G$, then it necessarily generates a subgroup isomorphic to $S_4$, 
but since elements of order $4$ in $S_4$ are odd
permutations while elements of order $3$ are even, this is a contradiction. 
Hence all such triples necessarily generate $G$.
Similarly, a non-generating triple associated with $(0;4,4,4)$
generates a subgroup isomorphic to a subgroup of $S_4$ or $D_8$,
but since elements of order $4$ in $S_4$ are odd, 
and $D_8$ has a unique cyclic subgroup of order $4$, this is
a contradiction as well, implying that triples of this shape are
necessarily generating.
 
\bigskip

\noindent
It remains to look at $(0; 3, 4, 4)$: The relevant class 
multiplication coefficient equals $c_{3, 4, 4} = 672$.
If an associated triple does not generate $G$, then it
necessarily generates a subgroup isomorphic to $S_4$.
Now the class multiplication coefficient in $S_4$ equals
$c_{3, 4, 4} = 42$, and since there precisely $7$ subgroups
in PSL$(2, 7)$ isomorphic to $S_4$, we conclude that there
are generating triples of this shape.
\qed

\bigskip

\noindent
Elaborating on the techniques indicated above, and also
using {\sf GAP} \cite{gap}, we have done a few computations
towards an understanding of the minimum stable genus
and the bad solutions in the general case.
The explicit results are shown in Table \ref{tbl4};
recall again the exceptions PSL$(2, 2)\cong S_3$ and PSL$(2, 3)\cong A_4$,
as well as PSL$(2, 4)\cong$ PSL$(2, 5)\cong A_5$,
and that we are only considering genera $g\geq 2$ here.

\begin{table}[here]
\begin{center}
\begin{tabular}{|r||r|r|r|r|}
\hline
        $p$ & $N_G$ & $\mu(G)$ & $\sigma(g)$ & bad \\
\hline
\hline
        $2$ & $1$ & $2$ & $2$ & $0$ \\
        $4$ & $1$ & $3$ & $63$ & $0$ \\
        $8$ & $2$ & $7$ & $1,453$ & $0$  \\
        $16$ & $4$ & $205$ & $32,153$ & $1$ \\
        $32$ & $8$ & $1,241$ & $517,617$ & $0$ \\
        $64$ & $16$ & $11,761$ & $1,386,081$ & $12$ \\
\hline
        $3$ & $1$ & $3$ & $3$ & $0$ \\
        $9$ & $3$ & $16$ & $505$ & $1$ \\
        $27$ & $9$ & $118$ & $61,696$ & $0$ \\
        $81$ & $27$ & $15,499$ & $5,371,111$ & $6$ \\
\hline
        $5$ & $1$ & $3$ & $63$ & $0$ \\
        $25$ & $5$ & $326$ & $52,111$ & $3$ \\
        $125$ & $25$ & $11,626$ & $9,886,176$ & $4$ \\
\hline
        $7$ & $1$ & $3$ & $399$ & $0$ \\
        $49$ & $7$ & $2,451$ & $337,359$ & $7$ \\
\hline
%
\end{tabular}
\end{center}
\caption{}\label{tbl4}
\end{table}

\noindent
These data imply that a straightforward generalization of Theorem 4.5
cannot possibly hold, leading to the following question:

\bigskip

\noindent {\bf Question 6:}
Give a closed formula for the stable upper genus of PSL$(2, q)$, 
and find a combinatorial description of the bad solutions.

\bigskip

\noindent
Finally, similarly, and still employing {\sf GAP} \cite{gap}, 
we have calculated the stable upper genus and the bad solutions
for the first thirteen
sporadic simple groups. The result is given in Table \ref{tbl5},
together with the respective genus increments and number of bad solutions. 

\begin{table}[here]
\begin{center}
\begin{tabular}{| l || r | r | r | r |}
\hline
        $G$ & $N_G$ & $\mu(G)$ & $\sigma(g)$ & bad \\
\hline
\hline
        M$_{11}$ & $3$ & $631$ & $48,511$ & $2$ \\
\hline
        M$_{12}$ & $36$ & $3,169$ & $510,841$ & $4$ \\
\hline
	J$_1$ & $2$ & $2,091$ & $2,749,249$ & $0$ \\
\hline
	M$_{22}$ & $24$ & $34,849$ & $3,856,897$ & $1$ \\
\hline 
	J$_2$ & $360$ & $7,201$ & $1,905,841$ & $3$ \\
\hline
	M$_{23}$ & $24$ & $1,053,361$ & $176,488,081$ & $2$ \\
\hline 
	HS & $2,400$ & $1,680,001$ & $335,793,601$ & $3$ \\
\hline
	J$_3$ & $216$ & $1,255,825$ & $880,271,713$ & $1$ \\
\hline
	M$_{24}$ & $576$ & $10,200,961$ & $4,063,754,881$ & $0$ \\
\hline
        McL & $16,200$ & $78,586,201$ & $6,587,730,001$ & $1$ \\
\hline
        He & $141,120$ & $47,980,801$ & $36,015,517,441$ & $1$ \\
\hline
        Ru & $115,200$ & $1,737,216,001$ & $2,658,295,065,601$ & $0$ \\
\hline
        Suz & $622,080$ & $11,208,637,441$ & $5,213,968,496,641$ & $0$ \\
\hline
\end{tabular}
\end{center}
\caption{}\label{tbl5}
\end{table}

\noindent 
This leads to the following question:

\bigskip

\noindent {\bf Question 7:} 
Compute the stable upper genus and the number of bad solutions
for the remaining sporadic simple groups.

\vspace{.4in}


\vspace*{.4in}

\noindent
{\sc J.M.:
Lehrstuhl D f\"ur Mathematik, RWTH Aachen \\
Templergraben 64, D-52062 Aachen, Germany} \\
{\sf juergen.mueller@math.rwth-aachen.de}

\bigskip

\noindent
{\sc S.S.:
Mathematics Department, IISER Bhopal \\
Indore By-pass Road, Bhauri \\
Bhopal 462030 (M.P.), India} \\
{\sf sidhu@iiserb.ac.in}

\end{document}